\documentclass[10pt]{amsart}
\usepackage{amsmath}
\usepackage{amssymb}
\usepackage{amsthm}
\usepackage{amscd}
\usepackage{amsfonts}
\usepackage{graphicx}
\newcommand{\binombf}[2]{\left( \begin{matrix}
  {\mathbf #1} \\
  {\mathbf #2}
\end{matrix}\right)}
\newcommand{\La}{\mathcal{L}_\nu}
\newcommand{\LO}{L^2_{\nu}(\Omega)}

\newcommand{\TO}{\mathop{T(\Omega)} }

\newcommand{\RR}{\mathbb R}

\newcommand{\CC}{\mathbb C}
\newcommand{\fk}{\mathfrak{k}}
\newcommand{\fg}{\mathfrak{g}}
\newcommand{\fp}{\mathfrak{p}}

\newcommand{\SO}{\mathop{\mathrm{SO}}}
\newcommand{\SL}{\mathop{\mathrm{SL}}}

\newcommand{\Sp}{\mathop{\mathrm{Sp}}}
\newcommand{\GL}{\mathop{\mathrm{GL}}}

\newcommand{\lsl}{\mathop{\mathfrak{sl}}}
\newcommand{\lsp}{\mathop{\mathfrak{sp}}}
\newcommand{\hO}{\mathcal{H}_\nu (T (\Omega ))}

\newcommand{\abs}[1]{\left\vert#1\right\vert}
\newcommand{\set}[1]{\left\{#1\right\}}

\newcommand{\ip}[2]{\left(#1\,\vert\, #2 \right)}

\newcommand{\RE}{\mathrm{Re}}

\newcommand{\ad}{{\operatorname{ad}}}
\newcommand{\lspan}{{\operatorname{span}}}
\theoremstyle{plain}
\newtheorem{thm}{Theorem}[section]
\newtheorem{cor}[thm]{Corollary}
\newtheorem{lem}[thm]{Lemma}
\newtheorem{prop}[thm]{Proposition}

\theoremstyle{definition}

\newtheorem{rem}[thm]{Remark}

\def\tr{\operatorname{tr}}

\numberwithin{equation}{section}

\begin{document}
\title[Laguerre Functions on Symmetric Cones]
{Laguerre Functions on Symmetric Cones and recursion relations in the Real Case}
\author{Michael Aristidou, Mark Davidson, and Gestur \'Olafsson }
\address{Department of Mathematics, Louisiana State University, Baton Rouge, LA\ 70803, USA}
\email{aristido@math.lsu.edu}
\email{davidson@math.lsu.edu}
\email{olafsson@math.lsu.edu}
\thanks{Aristidou was partially supported
by NSF grant DMS-0139783.
Both
Aristidou and \'Olafsson were partially supported by the GSF - National Research Center for
Environment and Health.
The research of \'Olafsson was supported by NSF grants DMS-0139783 and DMS-0402068.}

\begin{abstract}
In this article we derive differential recursion relations for the
Laguerre functions on the cone $\Omega$ of positive definite real
matrices. The highest weight representations of the group $\Sp
(n,\RR)$ play a fundamental role. Each such representation  acts on
a Hilbert space of holomorphic functions on the tube domain $\Omega
+ i\mathrm{Sym}(n,\RR)$. We then use the Laplace transform to carry
the Lie algebra action  over to $L^2(\Omega ,d\mu_\nu)$. The
differential recursion relations result by restricting to a
distinguished three dimensional subalgebra, which is isomorphic to
$\mathfrak{sl}(2,{\mathbb R}).$
\end{abstract}

\maketitle

\section*{Introduction}
\noindent The theory of special functions has its origins in the
late eighteenth and early nineteenth centuries when it was seen that
the algebraic, exponential, and trigonometric functions (and their
inverses)  were not adequate to express results to differential
equations that arose in the context of some important physical
problems. New functions arose to which  we have associated names
like Bessel, Hermite, Jacobi, Laguerre, and Legendre. Then there are
the  Gamma, Beta, Hypergeometric and many other families of special
functions.   By the latter half of the $19^{\text{th}}$ century
these same functions arose in different contexts and their name
`special' began to take on greater meaning. Their functional
properties were explored and included functional relations,
differential and difference recursion relations, orthogonality
relations,  integral relations, and others.

The preface of Vilenkin's  book \cite{V68} notes that the
connection between special functions and group representations was
first discovered by \`E. Cartan in the early part of the
$20^{\text{th}}$ century.  By the time Vilenkin's book appeared in
the 1960's that interplay had been well established.  The texts by
Miller, Vilenkin, and  Vilenkin and Klimyk \cite{M,V68,VK91}, for
example, well document the general philosophy. In short, group
representation theory   made it possible to express the classical
special functions as matrix entries of a representation and  to
unify many of the disparate relationships mentioned above. The
representation can then be used to derive differential equations
and differential recursion relations for those functions.

In 1964 Simon Gindikin published his paper: `Analysis on Homogeneous
Domains', cf. \cite{G64}.  This important paper developed special
functions as part of analysis on homogeneous convex cones and built
upon the earlier work of C. L. Siegel \cite{Siegel} on the cone of
positive definite matrices. The Siegel integral of the first and
second kind generalize to become the Beta and Gamma functions for
the cone, respectively. Generalized hypergeometric functions  are
extended to homogeneous cones. Many important differential
properties also extend.

Around the same time M. Koecher \cite{koe58a,koe62b}  began to
develop his analysis on symmetric cones and the complex tube domains
associated with them. Jordan algebras proved to be a decisive tool
for framing and obtaining many important fundamental results.  The
outstanding text by Faraut and Koranyi \cite{fk} documents this
interaction (see also its extensive bibliography). Nevertheless, the
representation theory of Hermitian  groups, which are naturally
associated with tube domains, is not used in any outstanding way.

In a series of papers \cite{doz1,doz2,do,do04} the second and third
authors (with Genkai Zhang in the first two referenced articles) use
the representation theory of Hermitian groups in a decisive way to
obtain differential and difference recursion relations on series of
special functions. In the context of bounded symmetric domains the
relevant special functions are generalized Meixner polynomials and
in the context of tube domains over a symmetric cone the relevant
special functions are Laguerre functions. These special functions
exist in distinguished $L^2$-spaces, which are unitarily isomorphic
to Hilbert spaces of holomorphic functions on either a bounded
symmetric domain or a tube domain. The well-known representation
theory that exists there then transfers to the corresponding
$L^2$-space to produce differential and difference relations that
exist among the special functions. One cannot downplay the essential
role that Jordan algebras play in establishing and expressing many
of the fundamental results obtained about orthogonal families of
special functions defined on symmetric cones. Nevertheless, the
theory of highest weight representations add fundamental new results
not otherwise easily obtained.

In this present paper we will continue the themes outlined in the
above mentioned  papers for the Laguerre functions defined of the
cone of positive definite real symmetric matrices. The underlying
group is $\Sp(n,{\mathbb R})$ and its representation theory
establishes new differential recursion relations that Laguerre
functions satisfy. The case $n=1$  reduces to the classical Laguerre
functions defined on ${\mathbb R}^+$.   Briefly, the classical
Laguerre polynomials are defined by the formula $$L_m^\nu(x)=
 \sum_{k=0}^{m} \frac{\Gamma(m+\nu)}{\Gamma(k+\nu)}
\binom{m}{k} (-x)^k$$ and the Laguerre functions are defined by
$$\ell_m^\nu(x)=e^{-x}L_m^\nu(2x).$$ The classical differential
recursion relations are then expressed by the following three
formulas:\medskip
\begin{enumerate}
  \item $xD^2 + \nu D -x)\ell_m^\nu(x)=-(2m+\nu)\ell_m^\nu(x)$
  \item $xD^2 + (2x+\nu)D+
  (t+\nu))\ell_m^\nu(x)=-2m(\nu+m-1)\ell_{m-1}^\nu(x)$
  \item $xD^2 - (2x-\nu)D+
  (t-\nu))\ell_m^\nu(x)=-2\ell_{m+1}^\nu(x).$
\end{enumerate}\medskip
It is these three formulas that we generalize via the representation
theory of $Sp(n,{\mathbb R})$  to Laguerre functions defined on the
cone of positive definite real symmetric matrices. (Some definitions
of $L_m^\nu$ include a factor of $1/n!$. This is the case in
\cite{doz1} but not \cite{doz2}. The inclusion of this factor
changes the differential recursion relations slightly.)

This article is organized as follows: In the first section we
introduce some standard Jordan algebra notation. In particular, we
introduce the Laguerre functions and polynomials. Even though this
material and  most of the material in  Sections 2 and 3 hold in
general for simple Euclidean Jordan algebras, we specialize to the
case $J=\mathrm{Sym}(n,\RR )$, the Jordan algebra of symmetric real
$n\times n$-matrices. In Section 2 we introduce the tube domain
$\TO=\Omega +i\mathrm{Sym}(n,\RR )$, where $\Omega$ is the open self
dual cone of positive definite matrices. We also discuss the
structure of the group $\Sp (n,\RR )$ and its Lie algebra $\mathfrak
{sp}(n,{\mathbb R})$. Some important subalgebras of $\lsp (n, \CC)$
are introduced. This structure is later used to construct the
differential operators that give rise to the differential equations
satisfied by the Laguerre functions.

Section 3 is devoted to the discussion of the highest weight
representations $(\pi_\nu ,\hO)$. We also introduce the Laplace
transform as a special case of the \textit{restriction principle}
introduced in \cite{OOe96}. In Section 4 we  describe how the Lie
algebra, $\mathfrak{sp}(n,{\mathbb C})$ acts on $\hO$. In
particular, Proposition \ref{gactionprop} gives an  explicit formula
for the action of the derived representation for each of the three
subalgebras ${\mathfrak k}_{\mathbb C}$, ${\mathfrak p}^+$, and
${\mathfrak p}^-$, whose direct sum is $\mathfrak{sp}(n,{\mathbb
C})$. It should be noted, however, that not all of this information
is needed  to establish the differential recursion relations for the
Laguerre functions. In fact, only the action of three elements are
needed. The action of $\lsp (n,\CC )$ on $L^2(\Omega ,d\mu_\nu)$ is
described in Section \ref{ss:L2realization}. The  result is the
following theorem:
\medskip

\noindent \textbf{Theorem} \ref{algactthm}. \textit{ For $f\in
L^2_{\nu}(\Omega)$ a smooth vector we have:
\begin{enumerate}
\item$\lambda_{\nu}(X)f(x)=\tr[(bx+(ax-xa-\nu b)\nabla-x\nabla b\nabla]f(x)$, $X=\begin{pmatrix}
a & b\\
b & a
\end{pmatrix}\in\mathfrak{k}_\mathbb{C}$
\item$\lambda_{\nu}(X)f(x)=\tr[(\nu a+ax+(ax+xa+\nu a)\nabla+x\nabla a\nabla]f(x)$, $X=\begin{pmatrix}
a &  a\\
-a & -a
\end{pmatrix}\in\mathfrak{p}^+$
\item$\lambda_{\nu}(X)f(x)=\tr[(\nu a-ax+(ax+xa-\nu a)\nabla-x\nabla a\nabla]f(x)$, $X=\begin{pmatrix}
a & -a\\
a &-a
\end{pmatrix}\in\mathfrak{p}^-$
\end{enumerate}
}
\medskip

Here we use $\nabla$ to denote the gradient, $\mathfrak{k}_\CC$ is
the complexification of the Lie subalgebra $\mathfrak{u}(n)\subset
\lsp (n,\RR )$, and $\mathfrak{p}^\pm$ are certain Abelian
subalgebras of $\lsp (n,\CC )$ on which $\fk_\CC$ acts. We explain
the main ideas for the special case of $\Sp (1,\RR)\simeq \SL (2,\RR
)$ and present the lengthy proof of the theorem in the Appendix.

Specializing the above results to the elements
$$\textsc{x}=\begin{pmatrix} 1 & 1\cr -1 & -1\end{pmatrix}\in {\mathfrak p}^+\, ,
\quad \textsc{y}=\begin{pmatrix} 1 & -1\cr 1 &
-1\end{pmatrix}\in{\mathfrak p}^-\, , \quad\text{and}\quad
\textsc{z}=\begin{pmatrix} 0 & 1\cr 1 & 0\end{pmatrix}\in{\mathfrak
k}_{\mathbb C}\, ,
$$ where $1$ stands for the $n\times n$ identity matrix,
using properties of highest weight representations, and employing
Lemma 5.5 of \cite{doz2} we get our main result:
\medskip

\noindent
\textbf{Theorem} \ref{recursionthm}.
\textit{
The Laguerre functions are related by the following differential
recursion relations:
}
\begin{enumerate}
\item $\mathrm{tr}(-x\nabla\nabla-\nu
\nabla+x)\ell_{\mathbf{m}}^{\nu}(x)=
(n\nu+2|\mathbf{m}|)\ell_{\mathbf{m}}^{\nu}(x).$\vskip.1in
\item $ \mathrm{tr}(x\nabla\nabla+(\nu I+2x)\nabla+(\nu I+x)
)\ell_{\mathbf{m}}^{\nu}(x)= -2\sum_{j=1}^{r}%
\begin{pmatrix}
\mathbf{m}\\
\mathbf{m-\mathbf{e}_{j}}%
\end{pmatrix}\
(m_{j}-1+\nu-(j-1))\ell_{\mathbf{m-\mathbf{e}_{j}}}^{\nu}(x)
$\vskip.1in
\item $ \mathrm{tr}(-x\nabla\nabla+(-\nu I+2x)\nabla+(\nu I-x)
)\ell_{\mathbf{m}}^{\nu}(x)=
2\sum_{j=1}^{r}c_{\mathbf{m}}(j)\ell_{\mathbf{m+\mathbf{e}
_{j}}}^{\nu}(x).$
\end{enumerate}
\medskip

\section{The Jordan Algebra of Real Symmetric Matrices}
\noindent In this section we introduce the Jordan algebra $J$ of
real symmetric matrices. We then discuss the space of
$L$-invariant polynomial functions on $J$ and the
$\Gamma$-function associated to the cone of symmetric positive
matrices. Finally we introduce the generalized Laguerre functions
and polynomials.

\subsection{The Jordan Algebra $J=\mathrm{Sym} (n,\RR )$}
We denote by $J$ the  vector space of all real symmetric $n\times n$
matrices. The multiplication $x\circ y= \frac{1}{2}(xy + yx)$ and
the inner product $\ip{x}{y}=\tr(xy)$  turn $J$ into a real
Euclidean  simple Jordan algebra. The determinant and trace
functions for $J$ are the usual determinant and trace of an $n\times
n$ matrix and will be denoted $\det$ and $\tr$, respectively.
Observe that $\dim J :=d = \frac{n(n+1)}{2}.$ Let $\Omega$ denote
the interior of the cone of squares: $\set{x^2\mid x\in J}$. Then
$\Omega$ is the set of all positive definite matrices in $J$. Let
$$H(\Omega)=\set{g\in \GL(J) \mid g\Omega=\Omega}$$
and let
$H$ be the connected component of the identity of $H(\Omega)$. Then
$H$ can be identified with $\GL(n,{\mathbb R})_+$ (where $+$
indicates positive determinant)  acting on $\Omega$ by the formula
$$g\cdot x  = gxg^t,\quad g\in H,\; x\in \Omega.$$ This action is
transitive and, since $\Omega$ is self-dual, it follows that $\Omega$
is a symmetric cone. Let $L$ be the stability subgroup of the
identity $e\in \Omega$. Then $L=\SO(n,{\mathbb R})$ and
\begin{equation}\label{eq:Omega}
\Omega\simeq
H/L\, .
\end{equation}

Let $E_{i,i}$ be the diagonal $n\times n$ matrix with $1$ in the
$(i,i)$-position and zeros elsewhere.  Then
$(E_{1,1},\ldots,E_{n,n})$ is a Jordan frame for $J$.  Let $J^{(k)}$
be the $+1$-eigenspace of the idempotent $E_{1,1}+ \cdots + E_{k,k}$
acting on $J$ by multiplication.  Each $J^{(k)}$ is a Jordan
subalgebra and we have
$$J^{(1)}\subset J^{(2)} \subset\cdots \subset J^{(n)}=J.$$ If $\det_k$ is the
determinant function for $J_k$ and $P_k$ is orthogonal projection of
$J$ onto $J^{(k)}$ then the function $\Delta_k(x)=\det_k P_k(x)$ is
the usual $k^{\text{th}}$ principal minor for an $n\times n$
symmetric matrix; it is homogeneous of degree $k$. In particular
$\Delta (x):=\Delta_n(x)=\det (x)$. Note also that
\begin{equation}\label{eq:Hdelta}
\Delta (h\cdot x)=\det (h)^2\Delta (x)\, ,\qquad \forall h\in H\, .
\end{equation}

For  ${\mathbf m}=(m_1,\ldots, m_n)\in {\mathbb C}^n$ we write
${\mathbf m}\ge 0$ if each $m_i$ is a nonnegative integer and
$m_1\ge m_2 \ge \cdots \ge m_n\ge 0$. We let
$\Lambda=\set{{\mathbf m}\mid {\mathbf m}\ge 0}$. For each
${\mathbf m}\in \Lambda$  define
$$\Delta_{\mathbf m}= \Delta_1^{m_1-m_2}\Delta_2^{m_2-m_3}\cdots
\Delta_{n-1}^{m_{n-1}-m_n}\Delta_n^{m_n}.$$ These are the
\textit{generalized power functions}. It is not hard to see that
the degree of $\Delta_{\mathbf m}$ is $\abs{{\mathbf
m}}:=m_1+\cdots +m_n$. Observe that each generalized power
function extends to a holomorphic polynomial on $J_{\mathbb
C}=\mathrm{Sym}(n,{\mathbb C})$ in a unique way.
\subsection{The $L$-invariant Polynomials}
For each  ${\mathbf m}\in \Lambda$ we define an
$L$-invariant polynomial, $\psi_{\mathbf m}$, on $J_{\mathbb C} $ by
$$ \psi_{\mathbf m}(z)=\int_{L}^{ }\Delta_{\mathbf m}(lz) \, dl,\quad
z\in J_{\mathbb C} $$ where $dl$ is normalized Haar measure on $L$.
A well known theorem of Schmid (c.f. \cite{Schmid}) gives
\begin{lem}
If $\mathcal{P}(J_{\mathbb C})$ is the space of all polynomial
functions on $J_{\mathbb C}$ and $\mathcal{P}(J_{\mathbb C})^L$
denotes the space of $L$-invariant polynomials then
$\set{\psi_{\mathbf m}\mid {\mathbf m}\ge 0}$ is a basis of
$\mathcal{P}(J_{\mathbb C})^L$. Furthermore,  if
$\mathcal{P}_k(J_{\mathbb C})^L$ denotes the space of $L$-invariant
polynomials of degree less than or equal $k$ then
$\set{\psi_{\mathbf m}\mid \abs{{\mathbf m}}\le k}$ is a basis of
$\mathcal{P}_k(J_{\mathbb C})^L$.
\end{lem}

This lemma implies among other things that $\psi_{\mathbf m}(e+x)$
is a linear combination of $\psi_{\mathbf n}$, $\abs{{\mathbf
n}}\le \abs{{\mathbf m}}.$  The \textit{generalized binomial
coefficients}, $\binombf{m}{n}$, are  defined by the
equation:$$\psi_{\mathbf m}(e+x)=\sum_{\abs{{\mathbf n}}\le
\abs{{\mathbf m}}}^{ }\binombf{m}{n}\psi_{\mathbf n}(x).$$

\subsection{The Generalized Gamma Function} For ${\mathbf m}\in {\mathbb
C}^n$ we define $\Delta_{\mathbf m}(x)$, $x\in\Omega$, by the same
formula given above for ${\mathbf m}\in \Lambda$.  The
\textit{generalized Gamma function} is defined by
$$ \Gamma_\Omega({\mathbf m})=\int_{\Omega}^{ }e^{-\tr
x}\Delta_{\mathbf m}(x) \Delta(x)^{\frac{-(n+1)}{2} }\, dx.$$
Conditions for convergence of this integral are given in the
proposition below. If $\lambda$ is a real number we will associate
the multi-index $(\lambda,\cdots,\lambda)$ and denote it by
$\lambda$ as well. The context of use should not cause confusion.
Thus we define
$$(\lambda)_{\mathbf m}=\frac{\Gamma_\Omega(\lambda + {\mathbf
m})}{\Gamma_\Omega(\lambda)}.$$

For later reference we note the following facts about the
generalized Gamma function:
\begin{prop}\label{gammaprop} Let the notation be as above. Then
the following holds:
\begin{enumerate}
  \item If ${\mathbf m}=(m_1,m_2,\ldots,m_n)\in \CC^n$ then the integral defining the
generalized  Gamma function converges if
$\RE(m_j)>(j-1)\frac{1}{2}$, for $j=1,\ldots ,n$, and in this case
$$\Gamma_\Omega({\mathbf
  m})=  (2\pi)^\frac{n(n-1)}{4} \prod_{i=1}^n \Gamma(m_i
  -(i-1)\frac{1}{2}),$$ where $\Gamma$ is the usual Gamma function.
  In particular it follows, that the $\Gamma$-function has
  a meromorphic continuation to all of $\CC^n$.
  \item If ${\mathbf e}_k$ is an $n$-vector with $1$ in the $k^{\text{th}}$
  position and $0$'s elsewhere
  then the following holds for all $\mathbf{m}\in\CC^n$:
  \begin{enumerate}
  \item  $\displaystyle{\frac{\Gamma_\Omega({\mathbf
  m})}{\Gamma_\Omega({\mathbf
  m}-\mathbf{e}_k)}= m_k-1   -(k-1)\frac{1}{2} }$;
\smallskip

  \item $\displaystyle{\frac{\Gamma_\Omega({\mathbf
  m}+\mathbf{e}_k)}{\Gamma_\Omega({\mathbf
  m})}= m_k   -(k-1)\frac{1}{2} }$.
\end{enumerate}
\end{enumerate}
\end{prop}
\begin{proof}
Part (2) follows immediately from (1) and part (1) is Theorem 7.1.1
of \cite{fk}.
\end{proof}

\subsection{The Generalized Laguerre Functions and Polynomials}
Let $\nu> 0$ and ${\mathbf m}\in \Lambda$. The \textit{generalized
Laguerre polynomial} is defined (cf. \cite{fk} p. 242) by the
formula
$$L^\nu_{\mathbf m}(x)=(\nu)_{\mathbf m} \sum_{\abs{{\mathbf n}}\le\abs{{\mathbf m}} }^{
} \binombf{m}{n} \frac{1}{(\nu)_{\mathbf n} }{\psi_{\mathbf n}(-x)}, \quad x\in J$$ and
the \textit{generalized Laguerre function} is defined by
$$\ell_{\mathbf m}^\nu(x)=e^{-\tr x}L_{\mathbf m}^\nu(2x).$$

\begin{rem}
In the case $n=1$, i.e. in the case $G\simeq Sp(1,{\mathbb
R})=\SL(2,{\mathbb R})$, the generalized Laguerre polynomials and
functions  defined above are precisely the classical Laguerre
polynomials and functions defined on ${\mathbb R}^+$. We refer to
\cite{doz1} for the discussion of that case.
\end{rem}
The determinant $\mathrm{Det}(h)$ of $h\in H$ acting on $J$ is
$$\mathrm{Det}(h)=\det (h)^{n+1}\, .$$
It follows from (\ref{eq:Hdelta}) that the measure
$$d\mu_0 (x)=\Delta (x)^{-\frac{n+1}{2}}\, dx$$
is $H$-invariant. Here $dx$ is the Lebesgue measure on $J$.
More generally, we set
$d\mu_\nu(x)=\Delta(x)^{\nu-\frac{n+1}{2}} dx$ and define
$$L^2_{\nu}(\Omega)=L^2(\Omega,d\mu_\nu).$$  We observe that by
(\ref{eq:Hdelta}) it follows that $H$ acts
unitarily on $L^2_{\nu}(\Omega)$ be the formula
$$\lambda_\nu(h)f(x)=\det(h)^\nu  f(h^t\cdot x).$$

\begin{thm}[\cite{doz2,fk}]
The set
$\set{\ell_{\mathbf m}^\nu\mid {\mathbf m}\ge 0}$
forms a complete orthogonal
system in $L^2_{\nu}(\Omega)^L$, the Hilbert space
of $L$-invariant functions in $L^2_{\nu}(\Omega).$
\end{thm}

In \cite{doz1} it was shown, that the classical differential
recursion relations and differential equations for the Laguerre
functions on $\RR^+$ follows from the representation theory of
$\Sp (1,\RR )= \SL (2,\RR )$. In \cite{do} this was generalized to
the space of complex Hermitian matrices. It is a goal of this
article to extend this result to the generalized Laguerre
functions defined of the cone of symmetric matrices. This
indicates nicely what the more general results should be. Here we
use heavily the structure of $\Sp (n,\RR )$ and its Lie algebra,
but the proof of the general results should be more in the line of
Jordan algebras.

\section{The Tube Domain $\TO$, the Group $\Sp (n,\RR)$
and its Lie Algebra}

\noindent In this section we introduce the tube domain $\TO=\Omega
+i\mathrm{Sym}(n,\RR )$ and discuss the action of the group $\Sp
(n,\RR )$ on this domain. We then discuss some important Lie
subalgebras of $\lsp (n,\CC )$, the complexification of the Lie
algebra of $\fg=\lsp (n,\RR )$. These subalgebras will show up again
in Section \ref{s:Lieaction} where we compute their action  on
Hilbert spaces of holomorphic functions on $\TO$ introduced in the
next section. We then use that information to construct the Laguerre
differential operators.

\subsection{The Group $\Sp (n,\RR )$}
Let $T(\Omega)=\Omega+ iJ$ be the tube over $\Omega$ in $J_{\mathbb C}$, which we
identify with the space of complex $n\times n$ symmetric matrices. Let $G_\circ$ be the
group of biholomorphic diffeomorphisms on $T(\Omega).$ Then $G_\circ$ is a Lie group
with Lie algebra isomorphic to $\mathfrak{sp}(n,{\mathbb R})$ and acts homogeneously on
$T(\Omega)$. The group $\Sp(n,{\mathbb R})$ is isomorphic to a finite covering group of
$G_\circ$ in the following precise way. The usual definition of $\Sp(n,{\mathbb R})$ is
$$\Sp(n,{\mathbb R})=\set{g\in \SL(2n,{\mathbb R})\mid gjg^t=j},$$ where $j=\begin{pmatrix}
  0 & 1 \\
  -1 & 0
\end{pmatrix}\in \SL(2n,{\mathbb R}).$  Each element $g\in \Sp(n,{\mathbb R})$ can be
written in block form as $$g=\begin{pmatrix}
  A & B \\
  C & D
\end{pmatrix},$$ where $A,B,C,$ and $D$ are real matrices.  Defined in this way $\Sp(n,{\mathbb
R})$ acts by linear fractional transformation on the \textit{upper
half plane} $J+i\Omega$.  Let $G$ be the group defined by
$$G=\set{\begin{pmatrix}
  A & B \\
  C & D
\end{pmatrix}\mid \begin{pmatrix}
  A & -iB \\
  iC & D
\end{pmatrix}\in \Sp(n,{\mathbb R})}.$$ This means then that the $(1,2)$ and $(2,1)$ entries of
an element of $G$ are purely imaginary matrices.  Clearly $G$ is
isomorphic to $\Sp(n,{\mathbb R})$.  It acts on the \textit{right
half plane} $T(\Omega)=\Omega+iJ$  by linear fractional
transformations: if $g=\begin{pmatrix}
  A & B \\
  C & D
\end{pmatrix}\in G$ and  $z\in T(\Omega)$ then $$g\cdot z =(Az+B)(Cz+D)^{-1}.$$ It is
also a finite covering group of $G_\circ$.  For $g=\begin{pmatrix}
  A & B \\
  C & D
\end{pmatrix}\in G$  we have the
following relations among $A,B,C$, and $D$:
$$\begin{matrix}
A^tC-C^tA=0 &\quad  &AB^t-BA^t=0 \\
A^tD-C^tB=I & \quad  &AD^t-BC^t=I \\
 B^tD-D^tB=0 & \quad  &  CD^t-DC^t=0 \\
 B^tC-D^tA=-I & \quad  &AD^t-BC^t=-I
\end{matrix}$$

\subsection{Some Subgroups of G} Let $e$ be the $n\times n$-identity matrix.
Then $e\in \Omega\subset T(\Omega)$. Let $K$ be the stability
subgroup of $e$ in $G$. Then
$$K=\set{\begin{pmatrix}
  A & B \\
  B & A
\end{pmatrix} \in G\mid A\pm B \in U(n) }\simeq U(n)\, .$$
Here the isomorphism is given by
$$
\begin{pmatrix}
  A & B \\
  B & A
\end{pmatrix} \mapsto A+B\, .$$
The subgroup  $K$ is a maximal
compact subgroup of $G$ and $G/K$ is naturally isomorphic to $T(\Omega)$ by the map
$gK\rightarrow g\cdot e$.  The connected component of the identity of the subgroup of
$G$ that leaves $\Omega$ invariant is isomorphic to $H$ via the map
$$h\rightarrow \begin{pmatrix}
  h & 0 \\
  0 & (h^t)^{-1}
\end{pmatrix}\, .$$
This map realizes $L$ as a subgroup of $G$ as well. In fact, we
have
$$L=H\cap K,$$ via the above isomorphism.

\subsection{Lie Algebras}\label{ss:liealg} If $P=\begin{pmatrix}
  1 & 0 \\
  0 & i
\end{pmatrix}$ then $G = P^{-1}\mathrm{\Sp}(n,\mathbb{R})P$.
From this it follows that the Lie algebra, ${\mathfrak g}$, of $G$
is given by $\mathfrak{g}=P^{-1}\mathfrak{sp}(n,\mathbb{R})P$, and
hence
\begin{equation}\label{eq:Lie algebra}
\mathfrak{g}=\left\{
\begin{pmatrix}
a & ib\\
-ic & -a^t
\end{pmatrix}\in \mathfrak{sl}(2n,\mathbb {C})
\mid a,b,c\,\,\text{real}, b=b^t,c=c^t\right\} .
\end{equation}
We define a \textit{Cartan involution} on $\fg$ by
$$
\theta(X)=-X^*\, . $$ It induces a decomposition of $\fg$ and
$\fg_\CC$ into $\pm 1$-eigenspaces, ${\mathfrak g}= {\mathfrak k}+
{\mathfrak p} $ and $\fg_\CC = \fk_{\CC} +\fp_{\CC}$. The
$+1$-eigenspace ${\mathfrak k}$ is the Lie algebra of $K$. These
spaces are given by
$$\mathfrak{k}=\left\{X\in\mathfrak{g}\mid\theta(X)=X \right\}=\left\{
\begin{pmatrix}
 a & ib\\
ib & a
\end{pmatrix}\in \mathfrak{sl}(2n,\mathbb {C})
\mid a,b\,\,\text{real}, a=-a^t,b=b^t,\tr(a)=0\right\}$$ and
$$\mathfrak{p}=\left\{X\in\mathfrak{g}\mid \theta(X)=-X \right\}=\left\{
\begin{pmatrix}
a &  ib\\
-ib & -a
\end{pmatrix}\in\mathfrak{sl}(2n,\mathbb {C})
\mid a,b\,\,\text{real}, a=a^t,b=b^t\right\}.$$

Their complexifications  are  given by
\begin{equation*}
\mathfrak{k}_{\CC}=\left\{
\begin{pmatrix}
a & b\\
b & a
\end{pmatrix}\in \mathfrak{sl}(2n,\mathbb {C})
\mid a,b\,\,\text{complex}, a=-a^t,b=b^t,\tr(a)=0\right\}
\end{equation*}
and
\begin{equation*}
\mathfrak{p}_{\CC}=\left\{
\begin{pmatrix}
a & b\\
-b & -a
\end{pmatrix}\in\mathfrak{sl}(2n,\mathbb{C})
\mid a,b\,\,\text{complex}, a=a^t,b=b^t\right\}.
\end{equation*}

It is clear that $K_\CC$ acts on $\fp_{\CC}$. This representation
decomposes into two parts. For that let  $\textsc{z}=\begin{pmatrix}
0 & 1\\
1 & 0
\end{pmatrix}$. Then $\textsc{z}\in\mathcal{Z}(\mathfrak{k}_\mathbb{C})$ and
$\ad(\textsc{z}):\mathfrak{g}_\mathbb{C}\rightarrow\mathfrak{g}_\mathbb{C}$
has eigenvalues $0, 2, -2$.  The  $0$-eigenspace is
$\mathfrak{k}_\mathbb{C}$, the $+2$-eigenspace is denoted by
$\mathfrak{p}^-$ and is given by $$\mathfrak{p}^-=\left\{
\begin{pmatrix}
a &  -a\\
a & -a
\end{pmatrix}\in \mathfrak{g}_\mathbb{C}\mid a=a^t \right\}\subset {\mathfrak p}_{\mathbb C},$$ and the  $-2$-eigenspace
is denoted by $\mathfrak{p}^+$ and is given by  $$\mathfrak{p}^+=\set{
\begin{pmatrix}
a &  a\\
-a & -a
\end{pmatrix}\in \mathfrak{g}_\mathbb{C}\mid a=a^t }\subset {\mathfrak p}_{\mathbb C}.$$
Each of the spaces $\fp^\pm$ are invariant under $K_{\CC}$ and
irreducible as $K_{\CC}$ representation.
Note that this is not necessarily the standard notation. In our
notation the eigenvectors (in ${\mathfrak p}^+$) with
$-2$-eigenvalue correspond to  \textit{annihilation} operators while
the eigenvectors (in ${\mathfrak p}^-$) with $+2$-eigenvalue
correspond to  \textit{creation} operators. These operators will be
described is Section \ref{s:Lieaction} below.

\section{The Highest Weight Representations $(\pi_\nu ,\hO )$}\label{s:hwrep}

\noindent
In this section we introduce the
highest weight representations
$(\pi_\nu ,\hO)$ and state the main results needed. We also introduce the
Laplace transform as a special case of the \textit{restriction principle}
introduced in \cite{OOe96}.
\subsection{Unitary Representations of $G$ in $\mathcal{O}(\TO )$}
In this subsection we define a series of unitary representations of
$G$ on a Hilbert space of holomorphic
functions on $T(\Omega)$.  These representations are well known.  Let $\tilde{G} $ be
the universal covering group of $G$. Then $\tilde{G}$ acts on $T(\Omega)$ by
$(g,z)\mapsto\kappa(g)\cdot z$ where $\kappa:\tilde{G}\rightarrow G$ is the canonical
projection. For $\nu> n$ let $\mathcal{H}_{\nu}(T(\Omega))$ be the space of holomorphic
functions $F:T(\Omega)\rightarrow\mathbb{C}$ such that
\begin{equation}
\|  F\|
_{\nu}^{2}:=\alpha_{\nu}\int_{T(\Omega
)}|F(x+iy)|^{2}\Delta (x)^{\nu - (n+1)}\,dxdy<\infty \label{eq-normub}
\end{equation}
where
\begin{equation*}\label{eq:alpha}
\alpha_{\nu}=\frac{2^{n\nu}}{(4\pi)^{d}\Gamma_{\Omega}(\nu-\frac{n+1}{2} )}\,.
\end{equation*}
Then $\mathcal{H}_{\nu}(T(\Omega))$ is a non-trivial Hilbert space
with inner product
\begin{equation}\label{eq:ipr}
\ip{F}{G}=\alpha_\nu\int_{\TO} F(x+iy)\overline{G(x+iy)}\,
\Delta (x)^{\nu - (n+1)}\, dxdy\, .
\end{equation}
For $\nu\le n$ this space reduces to $\set{0}$. If
$\nu=n+1$ this is the \textit{Bergman space}.  The space
$\mathcal{H}_\nu(T(\Omega))$ is a \textit{reproducing kernel
Hilbert space}. This means that point evaluation
$$E_z:\mathcal{H}_\nu((T(\Omega))\rightarrow {\mathbb C}$$ given
by $E_zF=F(z)$ is continuous, for every $z\in T(\Omega)$.  This
implies the existence of a kernel function $K_z\in \hO$, such that
$F(z)=\ip{F}{K_z}$ for all $F\in\hO$ and $z\in \TO$. Set $K(z,w)=K_w
(z)$. Then $K(z,w)$ is holomorphic in the first variable and
antiholomorphic in the second variable. The function $K(z,w)$ is
called \textit{the reproducing kernel} for $\hO$. We note that the
Hilbert space is completely determined by the function $K(z,w)$. In
particular, we have:
\begin{enumerate}
\item The space of finite linear combinations
$\hO^o:= \set{\sum c_j K_{w_j}\mid c_j\in\CC ,\, w_j\in \TO}$
is dense in $\hO$;
\item The inner product in $\hO^o$ is given by
$$\ip{\sum_j c_jK_{w_j}}{\sum_k d_k K_{z_k}}=
\sum_{j,k} c_j\overline{d_k}K (z_k ,w_j)\, .$$
\end{enumerate}
We refer to
\cite{doz2,Kunze} for more details.

\subsection{The Unitary Representations $(\pi_\nu ,\hO )$}
For $g\in\tilde{G}$ and $z\in T(\Omega)$, let $J(g,z)$ be the
\textit{complex} Jacobian determinant of the action of
$g\in\tilde{G}$ on $T(\Omega)$ at the point $z$. We will use the
same notation for elements $g\in G$. A straightforward calculation
gives

$$J(g,z)= \det(Cz+D)^{-n-1}, \quad g=\begin{pmatrix}
  A & B \\
  C & D
\end{pmatrix}\in G \text{ and } z\in T(\Omega).$$  We also have
the \textit{cocycle relation}:
\[
J(ab,z)=J(a,b\cdot z)J(b,z)
\]
for all $a,b\in\tilde{G}$ and $z\in T(\Omega)$. It is well known
that for $\nu>n $ the formula
\begin{equation}
\pi_{\nu}(g)f(z)=J(g^{-1},z)^{\frac{\nu}{n+1}}f(g^{-1}\cdot z)
=\det (A-z\,C )^{-\nu} f(g^{-1}\cdot z) \label{eq-actionub}
\end{equation}
defines a unitary irreducible representation of $\tilde{G}$. In
\cite{fk}, \cite{RV76}, and \cite{W79}  it was shown that this
unitary representation $(\pi_{\nu },\mathcal{H}_{\nu}(T(\Omega)))$
has an analytic continuation to the half-interval
$\nu>(n-1)\frac{1}{2}$. Here the representation $\pi_{\nu}$ is given
by the same formula (\ref{eq-actionub}) but the formula for the
\textit{norm} in (\ref{eq-normub}) is no longer  valid.
 There are also finitely many
equidistant values of $\nu$ that give rise to unitary
representations, but they will not be of concern to us here.

 In the following theorem we summarize what we have discussed and collect additional
information from \cite{doz2} and \cite{fk} (cf. p. 260, in
particular, Theorem XIII.1.1 and Proposition XIII.1.2).

\begin{thm}\label{th:mainsec2}
Let the notation be as above. Assume that for $\nu> n $ then the following hold:
\begin{enumerate}
\item  The space $\mathcal{H}_{\nu}(T(\Omega))$ is a reproducing Hilbert space.

\item  The reproducing kernel of $\mathcal{H}_{\nu}(T(\Omega))$ is given by
\[
K_{\nu}(z,w)=\Gamma_{\Omega}(\nu)\Delta\left(  z+\bar{w}\right)  ^{-\nu}%
\]

\item  If $\nu>\frac{1}{2}(n-1)$ then there exists a Hilbert space
$\mathcal{H}_{\nu}(T(\Omega))$ of holomorphic functions on
$T(\Omega)$ \ such that $K_{\nu}(z,w)$ defined in (2) is the
reproducing kernel of that Hilbert space. The representation
$\pi_\nu$ defines a unitary representation of $\tilde{G}$ on
$\mathcal{H}_{\nu}.$

\item  If $\nu>\frac{1}{2}(n-1)$ then the functions
\[
q_{\mathbf{m}}^\nu(z):=\Delta(z+e)^{-\nu}\psi_{\mathbf{m}}\left(
\frac
{z-e}{z+e}\right),  \,\qquad\mathbf{m}\in\mathbf{\Lambda},%
\]
form an orthogonal basis of $\mathcal{H}_{\nu}(T(\Omega))^{L}$, the space of
$L$-invariant functions in $\mathcal{H}_{\nu}(T(\Omega))$.
\end{enumerate}
\end{thm}

\subsection{The Restriction Principle and the Laplace Transform}
The restriction
principle \cite{OOe88,o00} is a general recipe to construct
unitary maps between a reproducing kernel Hilbert space
of holomorphic functions and $L^2$-spaces on a
totally real submanifold. Suppose $M$ is
a complex manifold and ${\mathbb H}(M)$ is a reproducing kernel
Hilbert space of holomorphic functions on $M$ with kernel $K$.
Suppose $X$ is a totally real submanifold of $M$ and a measure space
for some measure $\mu$.

Assume we have a holomorphic function
$D$ on $M$, such that $D$ is positive on $X$, and
such that the map
$$R:{\mathbb H}(M)\rightarrow L^2(X,\mu)\, ,$$
given by $Rf(x)=D(x)f(x)$, is densely defined.  As each $f$ is
holomorphic, its restriction to $X$ is injective. It follows that
$R$ is an injective map. We call $R$ a \textit{restriction map}.
Assume $R$ is  closed and has dense range. If $K(z,w)=K_w(z)$ is the
reproducing kernel for $\mathbb{H}(M)$, and $f\in L^2(X, d\mu )$,
then
\begin{eqnarray*}
R^*f(z) &=&\ip{R^*f}{K_z}\\
&=&\ip{f}{RK_z}\\
&=&\int_X f(x)D(x)K(z,x)\, d\mu (x)\, .
\end{eqnarray*}
In particular, if we set $\Psi(x,y)=D(y)D(x)K(y,x)$,
then $RR^*$ is given by
$$RR^* f(y)=\int_X f(x)\Psi (x,y)\, d\mu (x)$$
and thus is  an integral operator. Consider the polar decomposition
of the operator $R^*$.  We can write $$R^*=U \sqrt{RR^*},$$ where
$U$ is a unitary operator
$$U:L^2(X,\mu)\rightarrow {\mathbb H}(M).$$ The unitary map $U$ is
sometimes called the \textit{generalized Segal-Bargmann transform}.
In many applications of the restriction principle,  $M$  and $X$
will be homogeneous spaces with a   group $H$ acting on both.  When
the restriction map $R$ is $H$-intertwining so will the unitary
operator $U$.  This is exactly what happens in the situation at
hand. Here we can take $D=1$ and define
$R:\mathcal{H}_{\nu}(T(\Omega))\rightarrow L^2_{\nu}(\Omega)$ by
$$Rf(x)=f(x)\, . $$
Then  we obtain the following:
 \begin{thm}
The map $R$ is injective, densely defined and has dense range. The
unitary part, $U$, of the polar decomposition of $R^*$:
$R^*=U\sqrt{RR^*}$, is the Laplace transform given by
$$Uf(z)=\mathcal{L}_{\nu}f(z)=\int_{\Omega}e^{-(z\mid x)}f(x)\,d\mu_{\nu}(x)\, .$$
Furthermore,
$$\mathcal{L}_{\nu}(\lambda_\nu(h)f)=\pi_\nu(h)\mathcal{L}_{\nu}(f)$$ for all $h\in H$.
In particular, $\mathcal{L}_\nu$ induces an isomorphism
$\mathcal{L}_\nu : L_\nu^2 (\Omega )^L\to \hO^L$. Moreover
$$\mathcal{L}_{\nu}(\ell_{\mathbf{m}}^{\nu})=\Gamma_{\Omega}(\mathbf{m}%
+\nu)q_{\mathbf{m}}^{\nu}\, .
$$
\end{thm}

\begin{proof}
The first proof of this theorem was done for $\SL(2,{\mathbb R})$ in
\cite{doz1}. The general case is on pages 187-190 of \cite{doz2}.
\end{proof}

\begin{rem}
Rossi and Vergne \cite{RV76} obtained the unitarity of the Laplace
transform using a result of Nussbaum.
\end{rem}

The unitarity of the Laplace transform allows us to transfer the
representation, $\pi_\nu$, of $G$ on $\mathcal{H}_{\nu}(T(\Omega))$
to an equivalent representation of $G$ on $L^2_{\nu}(\Omega)$, which
extends $\lambda_\nu$ by the above theorem, . We will denote the
extension by $\lambda_\nu$ as well.  It is possible to describe
$\lambda_\nu$ on various subgroups of $G$ whose product is dense in
$G$. However, it is a difficult problem at best to describe a global
realization of $\lambda_\nu$ on all of $G$. However, part of the
point of this paper is to give a formula for the derived
representation of $\lambda_\nu$ on the Lie algebra of $G$ and its
complexification. It is  from the derived representation that new
differential recursion relations arise that relate the generalized
Laguerre functions.

\section{The Action  of $\fg_\CC$}\label{s:Lieaction}

\noindent In this section we introduce some subalgebras of $\lsp
(n,\CC )$, the complexification of the Lie algebra of $G$, and
explain how they act in the Hilbert space $\hO$.

\subsection{The Derived Representation on $\hO$}
Denote by $\hO^\infty$ the space of functions $F\in \hO$ such that
the map
$$\RR\ni t\mapsto \pi_\nu (\exp tX )F\in \hO$$
is smooth for all $X\in \fg=\lsp (n,\RR )$. If $f\in C_c^\infty (G)$, then
$\pi_\nu (f)F=\int_G f(g)\pi_\nu (g)F\, dg$ is in $\hO^\infty$ and it follows,
that $\hO^\infty $ is dense in $\hO$.
The Lie algebra representation, denoted also by $\pi_\nu$,  of ${\mathfrak g}$ on
$\mathcal{H}_{\nu}(T(\Omega))^\infty $ is given, by differentiation as follows:
\begin{eqnarray*}
\pi_{\nu}(X)F&=&\lim_{t\to 0}\frac{\pi_\nu (\exp tX )F - F}{t}\\
&=&\frac{d}{dt}\pi_{\nu}(\exp(tX))F|_{t=0}.
\end{eqnarray*}
Note that the limit is taken in the Hilbert space norm in $\hO$, but
it is easy to see that if $F\in \hO^\infty$, then in fact for $X\in
\fg$:
\begin{equation}\label{eq:derived}
\pi_\nu (X)F(z)= \frac{d}{dt}J(\exp(-tX),z
)^{\frac{\nu}{n+1}}F(\exp(-tX)\cdot z)|_{t=0}\,,
\end{equation} for all $z\in T(\Omega).$
We extend this by complex linearity to ${\mathfrak g}_{\mathbb
C}$.

Define $D_w$ by
$$D_w F(z)=\frac{d}{dt}F(z+tw)\vert_{t=0}=F^\prime (z)w$$
where $F^\prime$ denotes the derivative of $F$.

\begin{lem}\label{detlem}
Suppose $z, w$ are $n\times n$ matrices over ${\mathbb C}$ and $z$
is invertible. Then
$$D_w \det(z)^n=n\det(z)^n \tr(z^{-1}w).$$
\end{lem}
\begin{proof}
This follows from the chain rule and the fact that $$D_w\det (z)
=\frac{d}{dt}\det(z+tw)\vert_{t=0}=\det z
\frac{d}{dt}\det(1+tz^{-1}w)\vert_{t=0}=\det (z) \tr(z^{-1}w)\, .$$
\end{proof}

The following proposition expresses the relevant formulas on
${\mathfrak k}_{\mathbb C}$, ${\mathfrak p}^+$, and ${\mathfrak
p}^-$.  It's proof is a straightforward calculation using   Lemma
\ref{detlem}.

\begin{prop}\label{gactionprop}
For each piece of the Lie algebra
of $\mathfrak{g}_\mathbb{C}$ introduced in
Subsection \ref{ss:liealg}, we have:
\begin{enumerate}
\item$\pi_{\nu}(X)F(z)=\nu\tr(bz)F(z)+D_{za-az-b+zbz}F(z)$, $X=\begin{pmatrix}
a & b\\
b & a
\end{pmatrix}\in\mathfrak{k}_\mathbb{C}$
\item$\pi_{\nu}(X)F(z)=-\nu\tr(az+a)F(z)-D_{(za+az)+zaz+a}F(z)$, $X=\begin{pmatrix}
a &  a\\
-a & -a
\end{pmatrix}\in\mathfrak{p}^+$.
\item$\pi_{\nu}(X)F(z)=-\nu\tr(-az+a))F(z)+D_{-(za+az)+zaz +a}F(z)$, $X=\begin{pmatrix}
a & -a\\
a & -a
\end{pmatrix}\in\mathfrak{p}^-$
\end{enumerate}
\end{prop}

\begin{rem} We note that these equations are closely related to
the Jordan algebra structure of $J$. In particular, we have that
$za+az=2z\circ a$, where $\circ$ denotes the Jordan algebra product
$a\circ b=\frac{1}{2}(ab+ba)$. Furthermore $zaz =P(z)a$ where $P$
denotes the quadric representation (cf. \cite{fk}, p 32). Finally we
note, that $\tr (az)=\tr (z\circ a)$. These formulas therefore
  clearly indicate that the more general results are expressible in
terms of  Jordan algebraic constructs.
\end{rem}

\subsection{Highest Weight Representations}  The fact that $\pi_\nu$ is a highest weight
representation plays a decisive role in the recursion relations that we obtain. At this
point we explain what this notion means.

We assume $G$ is a Hermitian  group, which means that $G$ is simple
and the maximal compact subgroup $K$ has a one dimensional center.
The   Hermitian groups have been classified in terms of their Lie
algebras. They are $\mathfrak{su}(p,q)$, $\mathfrak{sp}(n,{\mathbb
R})$, $ \mathfrak{so}^*(2n)$, $\mathfrak{so}(2,n)$, and two
exceptional Lie algebras. The assumption that $K$ has a one
dimensional center implies that $G/K$ is a bounded symmetric domain.
In particular $G/K$ is complex. It also implies  that the
complexification of the Lie algebra, ${\mathfrak g}_{\mathbb C}$,
has a decomposition of the form ${\mathfrak g}_{\mathbb
C}={\mathfrak p}^+\oplus {\mathfrak k}_{\mathbb C}\oplus {\mathfrak
p}^-$, Specifically, ${\mathfrak p}^+$, ${\mathfrak k}_{\mathbb C}$,
and ${\mathfrak p}^-$
 are the $-2$, $0$, $2$-eigenspaces of
$\ad(\textsc{z})$, respectively, where $\textsc{z}$ is in the center
of ${\mathfrak k}_{\mathbb C}$.

\begin{lem}\label{ppmlem}
We have the following inclusions:
\begin{eqnarray*}
&&[{\mathfrak k}_{\mathbb C},{\mathfrak p}^\pm]\subset {\mathfrak p}^\pm\\
&&[{\mathfrak p}^+,{\mathfrak p}^-]\subset{\mathfrak k}_{\mathbb C} \\
\end{eqnarray*}
\end{lem}

\begin{proof}
Suppose $X\in {\mathfrak p}^+$, $Y\in {\mathfrak p}^-$, and $Z\in
{\mathfrak k}_{\mathbb C}$.  Then
$$\ad(\textsc{z})[X,Y]=[\ad(\textsc{z})X,Y]+[X,ad(\textsc{z})Y]=-2[X,Y]+2[X,Y]=0.$$ This
implies $[X,Y]\in {\mathfrak k}_{\mathbb C}.$  Similarly,
$$\ad(\textsc{z})[Z,X]=[\ad(\textsc{z})Z,X]+[Z,\ad(\textsc{z})X]=-2[Z,X].$$ This implies
$[Z,X]\in \mathfrak{p}^+.$  A similar argument shows that $[Z,Y]\in
{\mathfrak p}^-$.
\end{proof}

Suppose $\pi$ that is an irreducible
representation of $G$ on a Hilbert Space
${\mathbb H}$. We say $\pi$ is a \textit{highest weight
representation} if there is a nonzero vector $v\in {\mathbb H}$ such
that
$$\pi(X)v=0,$$ for all $X\in {\mathfrak p}^+$.  Let ${\mathbb
H}_\circ$ be the set of all such vectors.  The following theorem
is well known.

\begin{thm}\label{th:hwrep}
Suppose $\pi$ is an irreducible  unitary highest weight
representation of $G$ on ${\mathbb H}$ and ${\mathbb H}_\circ$ is
defined as above. Then $(\pi\vert_K, {\mathbb H}_\circ)$ is
irreducible. Furthermore, there is a scalar $\lambda$ such that
$$\pi(\textsc{z})v=\lambda v,$$ for all $v\in {\mathbb H}_\circ$. If
$${\mathbb H}_n=\set{v\in {\mathbb H}\mid \pi(\textsc{z})v=(\lambda + 2n)v},$$
then $${\mathbb H}= \bigoplus_{n\ge 0} {\mathbb H}_n\, .$$
Additionally,
\begin{eqnarray*}
&&\pi(Z): {\mathbb H}_n\longrightarrow {\mathbb H}_{n}, \quad Z\in {\mathfrak k}_{\mathbb C}\, ; \\
 &&\pi(X): {\mathbb H}_n\longrightarrow {\mathbb H}_{n-1}, \quad X\in {\mathfrak p}^+\, ;\\
&&\pi(Y): {\mathbb H}_n\longrightarrow {\mathbb H}_{n+1}, \quad
Y\in{\mathfrak p}^-\, ,
\end{eqnarray*}
where, in the case $n=0$, ${\mathbb H}_{-1}$ is understood to be
the $\set{0}$ space.
\end{thm}
\begin{proof} By Lemma \ref{ppmlem}, ${\mathbb H}_\circ$ is an
invariant $K$-space. Suppose ${\mathbb V}_\circ$ is a nonzero
invariant subspace of ${\mathbb H}_\circ$ and ${\mathbb W}_\circ$ is
its orthogonal complement in ${\mathbb H}_\circ$.  Define ${\mathbb
V}_n$ inductively as follows: $${\mathbb V}_n=\lspan\set{\pi(Y)v\mid
Y \in{\mathfrak p}^-, v\in {\mathbb V}_{n-1}}.$$ Let ${\mathbb V}=
\oplus {\mathbb V}_n$.  Define ${\mathbb W}_n$ in the same way as
${\mathbb V}_n$ and let ${\mathbb W}=\oplus {\mathbb W}_n$. Then, by
Lemma \ref{ppmlem}, ${\mathbb V}$  and ${\mathbb W}$ are invariant
${\mathfrak g}_{\mathbb C}$ subspaces of ${\mathbb H}$.  Since $\pi$
is unitary ${\mathbb V}$ and ${\mathbb W}$ are orthogonal. However,
since $\pi$ is irreducible and ${\mathbb V}$ is nonzero, it follows
that  ${\mathbb V}={\mathbb H}$ and hence ${\mathbb W}=0$. This
implies ${\mathbb W}_\circ=0$ and thus $\pi\vert_K$ is irreducible.
Since $\pi(\textsc{z})$ commutes with $\pi(K)$ Schur's lemma implies
that $\pi(\textsc{z})=\lambda$ on ${\mathbb H}_\circ$ for some
scalar $\lambda$. Since ${\mathbb V}_\circ={\mathbb H}_\circ$,
induction, Lemma \ref{ppmlem}, and irreducibility of $\pi$ implies
that ${\mathbb V}_n={\mathbb H}_n$. The remaining claims follow from
Lemma \ref{ppmlem}.
\end{proof}

\begin{rem}
The operators $\pi(X)$, $X\in{\mathfrak p}^+$, are called
\textit{annihilation operators} because, for $v$ in the algebraic
direct sum $\bigoplus {\mathbb H}_n$, sufficiently many applications
of $\pi(X)$ annihilates $v$. For $Y\in {\mathfrak p}^-$ the
operators $\pi(Y)$ are called \textit{creation operators}.
\end{rem}
\begin{rem}
A straightforward calculation gives $$\pi_\nu(X)q_{\mathbf
0}^\nu=0,$$ for all $X\in{\mathfrak p}^+$ and that $\hO_\circ
={\mathbb C}q_{\mathbf 0}^\nu.$  Thus $(\pi_\nu,\hO)$ is an
irreducible unitary highest weight representation of $G$ and by
unitary equivalence so is $(\lambda_\nu,L^2_\nu(\Omega))$.

\end{rem}
\section{The Realization of $\lambda_\nu$ Acting on $L^{2}(\Omega,d\mu_{\nu })$}\label{ss:L2realization}
\noindent In this section we determine  explicitly the  action of
$\fg_\CC$ on  $L^{2}(\Omega,d\mu_{\nu })$. More specifically, we
define $\lambda_\nu$ via the Laplace transform by the following
formula $$\lambda_\nu(X)=  \mathcal{L}_{\nu}^{-1} \pi_\nu(X)
\mathcal{L}_{\nu}$$ and will determine explicit formulas for
$\lambda_\nu(X)$, for $X\in \mathfrak{p}^+$, $X\in {\mathfrak
k}_{\mathbb C}$, and $X \in {\mathfrak p}^-$.

\subsection{Preliminaries}
Let $E_{ij}$ be the $n\times n$ matrix with a $1$ in the $(i,j)$
position and $0$'s elsewhere.  Define
$\tilde{E}_{i,j}=\frac{1}{2}(E_{i,j}+E_{j,i}).$  Then the
collection $\set{\tilde{E}_{i,j}\mid 1\le i\le j\le n}$ is a basis
of $J$ and $J_{\mathbb C }$, the real and complex symmetric
matrices. Furthermore, $\ip{\tilde{E}_{i
,j}}{\tilde{E}_{k,l}}=\frac{1}{2}(\delta_{jk}\delta_{il}+\delta_{jl}\delta_{ik})$,
which implies $\set{\tilde{E}_{i,j}\mid 1\le i\le j\le n}$ is an
orthogonal basis. Set $D_{i,j}=D_{\tilde{E}_{i,j}}\, $ and observe
that $D_{i,j}=D_{j,i}.$  The \textit{gradient} of $f$, $\nabla f$,
is defined by
$$\ip{\nabla f(x)}{u}=D_uf(x).$$

\begin{prop}
Suppose $f,g\in L^2(\Omega, d\mu_{\nu})$ are smooth and $f$
vanishes on the boundary of the cone $\Omega$. Let $1\le i,j \le
n$. Then
\begin{enumerate}
\item
$$ \int_{\Omega}D_{i,j} f(s) g(s) ds=-\int_{\Omega}f(s) D_{i,j}
g(s) ds. $$
\item
$$ \int_{\Omega} e^{-(z\vert s)}z_{i,j} f(s) ds = \int_{\Omega} e^{-(z\vert
s)}D_{i,j} f(s) ds. $$
\end{enumerate}
\end{prop}

\begin{proof}
(1) is Stokes Theorem and (2) follows from (1) and the fact that
$D_{i,j}e^{-(z\vert s)}=-e^{-(z\vert s)}z_{i,j}$, $z\in J_{\mathbb
C}$.
\end{proof}

\subsection{The Representation $\lambda_\nu$}
Recall that we determined the action of $\fk_\CC$, $\fp^+$ and
$\fp^-$ on $\hO^\infty$ in Proposition \ref{gactionprop}. We
denote the subspace of smooth vectors in $L^2_\nu (\Omega )$ by
$L^2_\nu (\Omega )^\infty$. Thus $f\in L^2_\nu (\Omega )^\infty $
if and only if the map
$$\RR \ni t\mapsto \lambda_\nu (\exp tX)f\in L^2_\nu (\Omega )$$
is smooth for all $X\in \fg$. Thus
$$L^2_\nu (\Omega )^\infty =\mathcal{L}_\nu^{-1}(
\hO^\infty)\, .$$ The action of $\fg$ on $L^2_\nu (\Omega )^\infty
$ is, as usual, defined by
$$\lambda_\nu (X)f=\lim_{t\to 0}\frac{\lambda_\nu (\exp
tX)f-f}{t},$$ for $X\in {\mathfrak g}$,  and then by complex
linearity the action extends to $\fg_\CC$. The following theorem
collects the corresponding equivalent action on the Hilbert space
$L^2_\nu (\Omega )^\infty$. We  remark again that these formulas can
be stated in terms of the Jordan algebra structure of $J$ indicating
the extension of these results to other tube domains.
\begin{thm}\label{algactthm}
For $f\in L^2_{\nu}(\Omega)$ a smooth function we have:
\begin{enumerate}
\item$\lambda_{\nu}(X)f(x)=\tr[(bx+(ax-xa-\nu b)\nabla-x\nabla b\nabla]f(x)$, $X=\begin{pmatrix}
a & b\\
b & a
\end{pmatrix}\in\mathfrak{k}_\mathbb{C}$
\item$\lambda_{\nu}(X)f(x)=\tr[(\nu a+ax+(ax+xa+\nu a)\nabla+x\nabla a\nabla]f(x)$, $X=\begin{pmatrix}
a &  a\\
-a & -a
\end{pmatrix}\in\mathfrak{p}^+$
\item$\lambda_{\nu}(X)f(x)=\tr[(\nu a-ax+(ax+xa-\nu a)\nabla-x\nabla a\nabla]f(x)$, $X=\begin{pmatrix}
a & -a\\
a &-a
\end{pmatrix}\in\mathfrak{p}^-$
\end{enumerate}
\end{thm}

\medskip
\subsection{The Case of $Sp(1,{\mathbb R})$}
The proof, which appears in the appendix, is very long and
computational. However, to convey the main ideas of the proof we
will discuss the simpler case of $\Sp(1,{\mathbb R})$, which is
isomorphic to $\SL(2,{\mathbb R})$. A detailed account of this case
(modelled on the upper half plane) is found in \cite{doz1}.

Let $G=\set{\begin{pmatrix}
  a & ib \\
  -ic & d
\end{pmatrix}\mid \begin{pmatrix}
  a & b \\
  c & d
\end{pmatrix}\in \SL(2,{\mathbb R})}\, .$  The group $G$
acts on the right half plane $T({\mathbb R}^+)$ by linear fractional
transformations. The complexification, ${\mathfrak g}_{\mathbb C}$,
of the Lie algebra of G is  ${\mathfrak sl}(2,{\mathbb C}),$ a three
dimensional Lie algebra spanned by
$$ \textsc{z}=\begin{pmatrix}
  0 & 1 \\
  1 & 0
\end{pmatrix},\quad \textsc{x}=\begin{pmatrix}
  1 & 1 \\
  -1 & -1
\end{pmatrix}, \quad \text{and}\quad \textsc{y}=\begin{pmatrix}
  1 & -1 \\
  1 & -1
\end{pmatrix}\, .$$
Proposition \ref{gactionprop} for this case reads:

\begin{thm}\label{sl2Raction}
The action of $\lsl (2,\RR )$ on the right half-plane is given by:
\begin{enumerate}
\item $\pi_\nu(\textsc{z}) F(z)=\nu z F(z) + (z^2-1)F'(z)$
\item $\pi_\nu(\textsc{x})F(z)= -\nu(z+1)F(z) - (z+1)^2F'(z)$
\item $\pi_\nu (\textsc{y})F(z)= \nu(z-1)F(z) + (z-1)^2F'(z)$
\end{enumerate}
\end{thm}
To find the corresponding action on $L^2_\nu ({\mathbb R }^+)$ we
must compute the operators that corresponds to $D_z$, $M_z$,
$M_{z^2}$, $M_z\circ D_z$ and $M_{z^2}\circ D_z$ in $L^2_\nu
({\mathbb R}^+ )^\infty$. Here $M$ stands for ``multiplication
operator''. To do this, requires several uses of integration by
parts, a special case of Stokes theorem. It was exactly this kind
of computation that was done in \cite{doz1} and we repeat it here:

For $D_z$ we have:
\begin{eqnarray*}
\frac{d\, }{dz}\mathcal{L}_\nu (f)(z)&=&
\int_0^\infty \frac{d e^{-zt} }{dz} f(t)t^{\nu -1}\, dt\\
&=&\mathcal{L}_\nu (-t f(t))
\end{eqnarray*}
Thus $D_z\longleftrightarrow M_{-t}$.

For $M_z$ we have
\begin{eqnarray*}
z\mathcal{L}_\nu (f)(z)&=&\int_0^\infty -\frac{de^{-zt}}{dt} f(t)t^{\nu -1}\, dt\\
&=&\int_0^\infty e^{-zt}\frac{d\, }{dt}(f(t)t^{\nu -1}) \, dt\\
&=&\int_0^\infty e^{-zt}(f^\prime (t)+\frac{\nu -1}{t}f(t))\,
t^{\nu -1}\, dt.
\end{eqnarray*}
Thus $M_z\longleftrightarrow D +M_{\frac{\nu-1}{t}}$.

We calculate $M_{z^2}$ similarly and get
$$M_{z^2}\longleftrightarrow D^2+\frac{2(\nu -1)}{t}D
+\frac{(\nu -1)(\nu -2)}{t^2} \, .$$ Thus:
\begin{lem}\label{sl2Raction2}
Let the notation be as above. Then the following holds:
\begin{enumerate}
\item $D_z\circ \La =-\La\circ (M_t)$;
\item $M_z\circ \La =\La\circ \left(D+\frac{\nu -1}{t}\right)$;
\item $M_z\circ D_z\circ \La =\La \circ   (-t D-\nu ) $;
\item $M_{z^2}\circ D_z\circ\La =\La \circ   \left(-t D^2-2\nu D-\frac{\nu (\nu -1)}{t}\right)$.
\end{enumerate}
\end{lem}
Combining Theorem \ref{sl2Raction} and Lemma \ref{sl2Raction2} gives
\begin{lem}\label{sl2ractionlem} Let the notation be as above. Then the
following holds:
\begin{enumerate}
\item $\lambda_\nu (\textsc{z})=-tD^2-\nu D +t$;
\item $\lambda_\nu (\textsc{x})= tD^2+(\nu +2t)D+(\nu +t)$;
\item $\lambda_\nu (\textsc{y})=-tD^2+(-\nu +2t)D+(\nu -t)$.
\end{enumerate}
\end{lem}
Note that this is Theorem \ref{algactthm} for this special case.

One more ingredient is necessary for determining the classical
recursion relations. This is a direct calculation and given in the
following lemma. We note at this point that such a direct
calculation is not done in the general case; deeper properties of
the representation theory must be used. (c.f. Proposition
\ref{qmprop}.)

\begin{lem} \label{qmlem} Let $q_m^\nu(z)
=(z+1)^{-\nu}\left(\frac{z-1}{z+1}\right)^m$. Then
\begin{enumerate}
  \item $\pi_\nu(\textsc{z})q_m^\nu= (\nu+2m)q_m^\nu$
  \item $\pi_\nu(\textsc{x})q_m^\nu= -2m q_{m-1}^\nu$
  \item $\pi_\nu(\textsc{y})q_m^\nu= 2(\nu+m) q_{m+1}^\nu$
\end{enumerate}

\end{lem}

 The combination of Lemmas \ref{sl2ractionlem} and \ref{qmlem} gives
 the classical recursion relations stated in the introduction and
  proves Theorem \ref{recursionthm} for the classical case.

\section{Differential Recursion Relations for $\ell_{m}^{\nu}$}
\noindent We now turn our attention to differential recursion
relations that exist among the generalized  Laguerre functions.
These relations are obtained by way of the highest weight
representation $\lambda_\nu$ and  generalize the classical case
mentioned in the introduction.

We begin with some preliminaries and a result found in \cite{doz2}.
First we notice that in general the Lie algebra $\fg_\CC$ does not
map $(\LO^\infty)^L$ into itself. For the Laguerre functions the
full Lie algebra is too big;  we will in fact only need the much
smaller Lie algebra $\fg_\CC^L$, which maps $(\LO^\infty)^L$ into
itself. It is well known, that in case $\fg$ is simple, then
$\fg_\CC^L\simeq \lsl (2,\CC )$. We choose
$\textsc{z},\textsc{x},\textsc{y}$ so that the isomorphism, which we
will denote by $\varphi$, is given by
$$\textsc{z}\mapsto \begin{pmatrix} 0 & 1\cr 1 & 0\end{pmatrix}\, ,\quad
\textsc{x}\mapsto \begin{pmatrix}
  1 & 1 \\
  -1 & -1
\end{pmatrix}, \quad \text{and}\quad \textsc{y}\mapsto \begin{pmatrix}
  1 & -1 \\
  1 & -1
\end{pmatrix}\, .$$
Furthermore, we can assume that $\varphi (X^t)= \varphi(X)^t$. This
shows that several calculations can in fact be reduced directly to
$\lsl (2,\CC )$. We will come back to that later.

Define $Z^0:=\frac{1}{2}(\textsc{x}+\textsc{y})$. Then $\varphi
(Z^0)=
\begin{pmatrix}
  1 & 0 \\
  0 & -1
\end{pmatrix}$
and $Z^0$ is in the center of ${\mathfrak h}$. For ${\mathbf m}\in
\Lambda$ let
$$c_{\mathbf m}(j)=\prod_{j\neq
k}\frac{m_j-m_k- \frac{1}{2}(j+1-k)}{m_j-m_k-\frac{1}{2}(j-k)}\, .$$
Then by Lemma 5.5 in \cite{doz2} we have:

\begin{prop}\label{qmprop}
The action of $\textsc{z}$ and $Z^0$ is given by:
\begin{enumerate}
\item $\pi_{\nu}(\textsc{z})q_{\mathbf{m}}^\nu=(n\nu+2\left|\mathbf{m}\right|
)q_{\mathbf{m}}^\nu$.\vskip.05in
\item $\pi_{\nu}(-2Z^{0})q_{\mathbf{m}}^\nu=\sum_{j=1}^{r}\binom{\mathbf{m}%
}{\mathbf{m}-\gamma_{j}}q_{{\mathbf{m}}-\mathbf{e}_{j}}^\nu-\sum_{j=1}^{r}%
(\nu+{m}_{j}-\frac{1}{2}(j-1)))c_{\mathbf{m}}(j)q_{{\mathbf{m}}%
+\mathbf{e}_{j}}^\nu$.
\end{enumerate}
\end{prop}

\begin{cor}\label{recrelcor}
Let the notation be as above. Then the following holds:
\begin{enumerate}
  \item $\lambda_\nu(\textsc{z})\ell_{\mathbf{m}}^\nu=(n\nu+2\vert
  \mathbf{m}\vert)\ell_{\mathbf{m}}^\nu.$\vskip.05in
  \item $\lambda_\nu(-2Z^\circ)\ell_{\mathbf{m}}^\nu = \sum
_{j=1}^{r}{\binom{\mathbf m}{{\mathbf m}-\mathbf{e}_{j}}}(m_{j}-1+\nu-(j-1))\ell_{\mathbf{m}-\mathbf{e}_{j}%
}^{\nu}-\sum_{j=1}^{r}c_{m}(j)\ell_{\mathbf{m}+\mathbf{e}_{j}}^{\nu}$.
\end{enumerate}
\end{cor}

\begin{proof}  This statement follows from Proposition \ref{qmprop}  and the
following facts:
 $\lambda_\nu(X)=\mathcal{L}_{\nu}^{-1}
\pi_\nu(X)\mathcal{L}_{\nu}$, $\mathcal{L}_{\nu}(\ell_{\mathbf
m}^\nu)= \Gamma_{\Omega}({\mathbf m }+ \nu) q_{\mathbf m}^\nu$, and
Proposition \ref{gammaprop}. In each of these formulas  if either
index ${\mathbf m}+\mathbf{e}_j$ or ${\mathbf m}- \mathbf{e}_j$ is
not in $\Lambda$ then it should be understood that  the
corresponding function  does not appear.
\end{proof}

\begin{thm}\label{recursionthm}
The Laguerre functions are related by the following differential
recursion relations:
\begin{enumerate}
\item $\mathrm{tr}(-x\nabla\nabla-\nu
\nabla+x)\ell_{\mathbf{m}}^{\nu}(x)=
(n\nu+2|\mathbf{m}|)\ell_{\mathbf{m}}^{\nu}(x)$.
\item $\mathrm{tr}(x\nabla\nabla+(\nu I+2x)\nabla+(\nu I+x)
)\ell_{\mathbf{m}}^{\nu}(x)= -2\sum_{j=1}^{r}%
\begin{pmatrix}
\mathbf{m}\\
\mathbf{m-\mathbf{e}_{j}}%
\end{pmatrix}\
(m_{j}-1+\nu-(j-1))\ell_{\mathbf{m-\mathbf{e}_{j}}}^{\nu}(x)
$\vskip.1in
\item $\mathrm{tr}(-x\nabla\nabla+(-\nu I+2x)\nabla+(\nu I-x)
)\ell_{\mathbf{m}}^{\nu}(x)=
2\sum_{j=1}^{r}c_{\mathbf{m}}(j)\ell_{\mathbf{m+\mathbf{e}
_{j}}}^{\nu}(x).$
\end{enumerate}
\end{thm}
\begin{proof}
If $\textsc{z}=\begin{pmatrix}
  0 & 1 \\
  1 & 0
\end{pmatrix}$ then substituting $a=0$ and $b=1$ into Proposition
\ref{algactthm}, part 1, gives
$$\lambda_\nu(\textsc{z})=\tr(-x\nabla\nabla -\nu\nabla +x).$$ Combining
this with part 1 of Corollary \ref{recrelcor} gives the first
formula.

Recall that $\textsc{x}\in \fp^+$ and $\textsc{y}\in \fp^-$.
According to Theorem \ref{th:hwrep} we have that $\lambda_\nu
(\textsc{x})\ell^\nu_{\mathbf{m}}$ has to be a linear combination of
$\ell^\nu_{\mathbf{m}^\prime}$, with $m^\prime_j\le m_j$ for all
$j$. Similarly, $\lambda_\nu (\textsc{y})\ell^\nu_{\mathbf{m}}$ has
to be linear combination of those $\ell^\nu_{\mathbf{m}^\prime}$
with $m^\prime_j\ge m_j$. The statement follows now from Corollary
\ref{recrelcor} and the fact that $2Z^0=\textsc{x} + \textsc{y}$.
\end{proof}

\section{Some Open Problems}
\noindent There are still several open question that require further
work. We mention three of these. One is a relation to the classical
Laguerre polynomials, the other two are natural generalization of
classical relations.

\subsection{Relation to Classical Laguerre Functions} Every
positive symmetric matrix $A$ can be written as $A=kDk^{-1}$, where
$k\in \SO (n)$ and $D=d(t_1,\ldots ,t_n)$ is a diagonal matrix with
$t_j>0$. Thus, if
$$\Omega_1=\set{d(\mathbf{t})\mid \mathbf{t}\in (\RR^+)^n}\simeq (\RR^+)^n$$
then
$$\Omega =L\cdot \Omega_1\, .$$
As the Laguerre functions are $L$-invariant, it follows that they
are uniquely determined by their restriction to $\Omega_1$. Let
$T(\Omega_1):=\set{d(\mathbf{x} )+id(\mathbf{y})\mid \mathbf{x}\in
(\RR^+)^n, \mathbf{y} \in \RR^n}$. Then $T(\Omega_1)\simeq (\RR^+ +
i\RR)^n$, and the group $\SL (2,\RR)^n$ acts transitively on the
right hand side. But it is well known, that  $\SL (2,\RR)^n$ can be
realized as a closed subgroup of $\Sp (n,\RR)$. It follows
therefore, that the generalized Laguerre functions can be written as
a finite linear combinations of products of classical Laguerre
functions. It is a natural problem to derive an exact formula.

\subsection{Relations in the $\lambda$-parameter}
It is well known that the classical Laguerre polynomials satisfy the
following relations:
\begin{eqnarray*}
xL^\lambda_n&=&(n+\lambda +1)L^{\lambda -1}_n-(n+1)L^{\lambda-1}_{n+1}\\
xL^\lambda_n&=& (n+\lambda) L^{\lambda -1}_{n-1} -(n-x)L^{\lambda-1}_n\\
xL^{\lambda-1}_n&=&L^\lambda_n-L^\lambda_{n-1}\, .
\end{eqnarray*}
In \cite{doz1} it was shown, that these relations follows directly
from the representation theory of $\mathfrak{sl}(2,\RR)$. It is
therefore natural to look for similar relations for the generalized
Laguerre polynomials and functions.

\subsection{Relations in the $x,y$ Parameters}
Several other classical relations should be extended to the general
case. We name here only the following
$$L^{\alpha + \beta +1}_m(x+y)=
\sum_{n=0}^m L_n^\alpha (x)L^\beta_{m-n}(y)\, .$$ This relation is
closely related to the decomposition of the tensor product of two
highest weight representations and we expect that a similar relation
can be derived also for the general case. Notice, however, that for
general Laguerre polynomials the right hand side is $L$-invariant in
the $x$ and $y$ variable while that is not the case on the left hand
side. Thus any generalization will involve a projection (averaging
over $L$) onto the $L$-invariant functions.

\newpage
\section{Appendix: Proof of Theorem (\ref{algactthm})}

\begin{proof}
We will prove the theorem for the case $x\in \mathfrak{p}^+$. The
other two cases are done similarly. For convenience we let
$m=\nu-\frac{n+1}{2}.$

 Let $x=\begin{pmatrix}
   a & -a \\
   a & -a \
 \end{pmatrix} \in\mathfrak{p}^+$.  By Proposition \ref{gactionprop}

$$\pi_{\nu}(x)\mathcal{L}_{\nu}f(z)=-{\nu}\tr(az+a)\mathcal{L}_{\nu}f(z)-
D_{az+a+zaz+za}\mathcal{L}_{\nu}f(z),\quad X\in\mathfrak{p}^+.
$$ Let

\begin{eqnarray*}
A&=&-D_{az+za}\mathcal{L}_{\nu}f(z)\\
B&=&-D_{a}\mathcal{L}_{\nu}f(z)\\
C&=&-D_{zaz}\mathcal{L}_{\nu}f(z)\\
\end{eqnarray*}

\subsection*{Calculation of A}
\begin{eqnarray*}
A&=&-D_{az+za}\mathcal{L}_{\nu}f(z)\\
&=&-\int_{\Omega}D_{az+za}e^{-(z|x)}f(x)\det(x)^m\;dx\\
&=&\int_{\Omega}e^{-(z|x)}(az+za|x)f(x)\det(x)^m\;dx\\
&=&\sum_{i,j,k}\int_{\Omega}e^{-(z|x)}(a_{ik}z_{kj}+z_{ik}a_{kj})(x_{ji}f(x)\det(x)^m) \;dx\\
&=&\sum_{i,j,k}\int_{\Omega}e^{-(z|x)}(a_{ik}D_{kj}+a_{kj}D_{ik})(x_{ji}f(x)\det(x)^m)\;dx\\
&=&\sum_{i,j,k}\int_{\Omega}e^{-(z|x)}(a_{ik}\frac{1}{2}(\delta_{kj}\delta_{ji}+\delta_{ki}\delta_{jj})
+a_{kj}\frac{1}{2}(\delta_{kj}\delta_{ii}+\delta_{ij}\delta_{ki})
)f(x)\det(x)^m\;dx\\
&&+\sum_{i,j,k}\int_{\Omega}
e^{-(z|x)}(a_{ik}x_{ji}D_{kj}f(x)+a_{kj}x_{ji}D_{ik}f(x))\det(x)^m\;dx\\
& &
+\sum_{i,j,k}\int_{\Omega}e^{-(z|x)}f(x)m\det(x)^{m}(a_{ik}x_{ji}\tr(x^{-1}\tilde{E}_{kj})+
a_{kj}x_{ji}\tr(x^{-1}\tilde{E}_{ik}))dx\\
&=&(n+1)\tr(a)\int_{\Omega}e^{-(z|x)}f(x)\det(x)^m\;dx\\
&&+ \int_{\Omega}
e^{-(z|x)}( (ax+xa)_{ik}D_{ik}f(x))\det(x)^m\;dx\\
&& +\int_{\Omega}e^{-(z|x)}f(x)m\det(x)^{m} \tr(x^{-1} (xa+ax))dx\\
&=&(2m + n+1)\tr(a)\int_{\Omega}e^{-(z|x)}f(x)\det(x)^m\;dx+
\int_{\Omega}
e^{-(z|x)}( \tr((ax+xa)\nabla f)(x))\det(x)^m\;dx\\
&=&2\nu\tr(a)\int_{\Omega}e^{-(z|x)}f(x)\det(x)^m\;dx+ \int_{\Omega}
e^{-(z|x)}( \tr((ax+xa)\nabla f)(x))\det(x)^m\;dx\\
&=&2\nu\tr(a)\mathcal{L}_{\nu}(f)(x) +
\mathcal{L}_{\nu}((tr(ax+xa)\nabla)f)(x)
\end{eqnarray*}

\subsection*{Calculation of B:}

\begin{eqnarray*}
B&=&-D_{a}\mathcal{L}_{\nu}f(z)\\
&=&\int_{\Omega}e^{-(z|x)}(a|x)f(x)\det(x)^mdx\\
&=&\int_{\Omega}e^{-(z|x)}\tr(ax)f(x)\det(x)^mdx\\
&=&\mathcal{L}_{\nu}(\tr(ax)f)(x)
\end{eqnarray*}

\subsection*{Calculation of C:}
\begin{eqnarray*}
C&=&-D_{zaz}\mathcal{L}_{\nu}f(z)\\
&=&\int_{\Omega}e^{-(z|x)}(zaz|x)f(x)\det(x)^mdx\\
&=&\sum_{i,j,k,l}\int_{\Omega}e^{-(z|x)}a_{kl}z_{ik}z_{lj}x_{ji}f(x)\det(x)^m dx\\
&=&\sum_{i,j,k}\int_{\Omega}e^{-(z|x)}a_{kl}z_{ik}D_{lj}(x_{ji}f(x)\det(x)^m)dx\\
&=&\sum_{i,j,k,l}\int_{\Omega}e^{-(z|x)}a_{kl}z_{ik}\frac{1}{2}(\delta_{lj}\delta_{ji}+
\delta_{li}\delta_{jj})f(x)\det(x)^mdx\\
&&+\sum_{i,j,k,l}\int_{\Omega}
e^{-(z|x)}a_{kl}z_{ik}x_{ji}(D_{lj}f)(x)\det(x)^mdx\\
&&+\sum_{i,j,k,l}\int_{\Omega}e^{-(z|x)}a_{kl}z_{ik}x_{ji}f(x)m\det(x)^{m}\tr(x^{-1}\tilde{E}_{lj})dx\\
&=&\sum_{i,k,l}\int_{\Omega}e^{-(z|x)}a_{kl}z_{ik}(\frac{n+1}{2}\delta_{li})f(x)\det(x)^mdx\\
&&+\sum_{i,j,k,l}\int_{\Omega}e^{-(z|x)}a_{kl}D_{ik}(x_{ji}D_{lj}f(x)\det(x)^m) dx\\
&&+\sum_{i,j,k,l}\int_{\Omega}e^{-(z|x)}a_{kl}z_{ik} x_{ji}f(x)m\det(x)^{m}\tr(x^{-1}\tilde{E}_{lj})dx\\
 &=& \frac{n+1}{2}\tr(az)\int_{\Omega}e^{-(z|x)}f(x)\det(x)^mdx\\
 &&+\sum_{i,j,k,l}\int_{\Omega}e^{-(z|x)}a_{kl}\frac{1}{2}(\delta_{ij}\delta_{ki}+
\delta_{ii}\delta_{kj})(D_{lj}f)(x)\det(x)^mdx\\
&&+\sum_{i,j,k,l}\int_{\Omega}e^{-(z|x)}a_{kl}x_{ji}(D_{ik}D_{lj}f)(x)\det(x)^mdx\\
&&+\sum_{i,j,k,l}\int_{\Omega}e^{-(z|x)}a_{kl}x_{ji}D_{lj}f(x)m\det(x)^{m}\tr(x^{-1}\tilde{E}_{ik})dx\\
&&+m\sum_{j,k}\int_{\Omega}e^{-(z|x)}a_{lk}z_{kl}\delta_{il}f(x)\det(x)^{m}dx\\
\end{eqnarray*}
\begin{eqnarray*}
&=& \frac{n+1}{2}\tr(az)\int_{\Omega}e^{-(z|x)}f(x)\det(x)^mdx\\
 &&+\sum_{i,j,k,l}\int_{\Omega}e^{-(z|x)}a_{kl}\frac{1}{2}(\delta_{ij}\delta_{ki}+
\delta_{ii}\delta_{kj})(D_{lj}f)(x)\det(x)^mdx\\
&&+\sum_{i,j,k,l}\int_{\Omega}e^{-(z|x)}a_{kl}x_{ji}(D_{ik}D_{lj}f)(x)\det(x)^mdx\\
&&+\sum_{i,j,k,l}\int_{\Omega}e^{-(z|x)}a_{kl}x_{ji}D_{lj}f(x)m\det(x)^{m}\tr(x^{-1}\tilde{E}_{ik})dx\\
&&+m\sum_{j,k}\int_{\Omega}e^{-(z|x)}a_{lk}z_{kl}\delta_{il}f(x)\det(x)^{m}dx\\
&=& \nu\tr(az)\int_{\Omega}e^{-(z|x)}f(x)\det(x)^mdx\\
&&+\nu \int_{\Omega}e^{-(z|x)}(\tr(a\nabla)f)(x)\det(x)^mdx\\
&&+ \int_{\Omega}e^{-(z|x)}((\tr(x\nabla a \nabla)f)(x)\det(x)^mdx\\
&=& \nu\tr(az)\mathcal{L}_{\nu}f(x) + \nu
\mathcal{L}_{\nu}(\tr(a\nabla)f)(x)+\mathcal{L}_{\nu}(\tr(x\nabla
a\nabla)f)(x)
\end{eqnarray*}

\begin{eqnarray*}
\pi_{\nu}(x)\mathcal{L}_{\nu}f(z)&=&-{\nu}\tr(az+a)\mathcal{L}_{\nu}f(z)-
D_{az+a+zaz+za}\mathcal{L}_{\nu}f(z),\quad X\in\mathfrak{p}^+\\
&=&-{\nu}\tr(az+a)\mathcal{L}_{\nu}f(z)+2\nu\tr(a)\mathcal{L}_{\nu}(f)(x)\\
&&+\mathcal{L}_{\nu}((tr(ax+xa)\nabla)f)(x)+\mathcal{L}_{\nu}(\tr(ax)f)(x)\\
&&+\nu\tr(az)\mathcal{L}_{\nu}f(x) +
\nu\mathcal{L}_{\nu}(\tr(a\nabla)f)(x)\\
&&+\mathcal{L}_{\nu}(\tr(x\nabla a\nabla)f)(x)\\
&=&\nu\tr(a)\mathcal{L}_{\nu}(f)(x)+\mathcal{L}_{\nu}((tr(ax+xa)\nabla)f)(x)\\
&&+\mathcal{L}_{\nu}(\tr(ax)f)(x) +\nu\mathcal{L}_{\nu}(\tr(a\nabla)f)(x)\\
&&+\mathcal{L}_{\nu}(\tr(x\nabla a\nabla)f)(x)\\
&=&\mathcal{L}_{\nu}(\tr(\nu a+ax+(ax+xa+\nu a)\nabla+x\nabla
a\nabla)f(x))
\end{eqnarray*}
Taking the inverse Laplace transform of each side gives the desired
result.

\end{proof}


\begin{thebibliography}{9}


\bibitem {doz1} M. Davidson, G. \'{O}lafsson, and Genkai Zhang, \emph{Laguerre
Polynomials, Restriction Principle, and Holomorphic Representations
of $\mathrm{SL}(2,\mathbb{R})$} Acta Applicandae Mathematicae
\textbf{71} (2002), 261--277.

\bibitem {doz2} M. Davidson, G. \'{O}lafsson, and Genkai Zhang,
\emph{Segal-Bargmann Transform on Hermitian Symmetric Spaces and
Orthogonal Polynomials}, J. Funct. Anal. \textbf{204} (2003),
157--195

\bibitem{do} M. Davidson, and G. \'{O}lafsson, \emph{Differential Recursion
Relations for Laguerre Functions on Hermitian Matrices} \, Integral
Transforms and Special Functions \textbf{14} (2003), No 4, 469--484.


\bibitem{do04} M. Davidson, and G. \'Olafsson:
\emph{The Generalized Segal-Bargmann  transform and Special
Functions} Acta Applicandae Mathematicae, \textbf{81} (2004), 29--50



\bibitem {fk}J. Faraut, and A. Kor{a}nyi, \emph{Analysis on Symmetric Cones},
Clarendon Press, 1994.


\bibitem{G64} S.G. Gindikin: \emph{Analysis on
homogeneous domains}, Uspekhi Math. Nauk \textbf{19} (1964) 3--92;
Russian Math. Surveys, \textbf{19} (4) 1--89


\bibitem{koe58a} M. Koecher: \emph{Analysis in reellen Jordan
Algebren}, Nach. Akad. Wiss. G\"{o}ttingen Math.-Phys., \textbf{K1
II a,} 67-74.

\bibitem{koe62b} M. Koecher: \emph{Jordon Algebras and their
Applications}. Lectrues notes, Univ. of Missesota, Minneapolis.


\bibitem{Kunze} R. Kunze, \textit{Positive
Definite Operator-Valued Kernels and Unitary Representations}, In:
Proceedings of the Conference on Functional Analysis at Irvine,
California, (1966), 235-247, Thompson Book Company.



\bibitem{M} W. Miller Jr.: \emph{Lie Theorey and Special Functions},
Academic Press, 1968.


\bibitem {o00}G. \'{O}lafsson, \emph{Analytic Continuation in Representation
Theory and Harmonic Analysis}. In: Ed. J. P. Bourguignon, T.
Branson, O. Hijazi, Global Analysis and Harmonic Analysis, Seminares
et Congres \textbf{4}, (2000), 201--2333. The French Math. Soc.


\bibitem{OOe88} G. \'Olafsson and B. {\O}rsted: \emph{The
holomorphic discrete series for affine symmetric spaces, I}, J.
Funct. Anal. \textbf{81} (1988), 126--159

\bibitem {OOe96}  G. \'Olafsson and B. {\O}rsted: \emph{Generalizations of
the Bargmann transform}, Lie theory and its applications in physics.
Proceedings of the international workshop, Clausthal, Gemany, August
14--17, 1995. (H.-D. Doebner et al, ed.), World scientific,
Singapore, 1996, pp. 3--14, MR. \textbf{99}c:22017


\bibitem {RV76} H. Rossi and M. Vergne: \emph{
Analytic continuation of the holomorphic discrete series of a
semisimple Lie group}, Acta Math.. \textbf{136} (1976), 1-59.


\bibitem {Schmid}W.~Schmid, \emph{Die {Randwerte} holomorpher Funktionen auf
hermitesch symmetrischen {R\"{a}{}{}umen}}, Invent. Math \textbf{9}
(1969), 61--80.

\bibitem{Siegel}C. L. Siegel: \emph{\"Uber die analytische Theorie
der quadratishen Formen}, Ann. of Math. \textbf{36} (1935), 527-606.

\bibitem{V68} N.J. Vilenkin: Special Functions and the Theory of
Group Representations. Translations of Mathematical Monographs,
\textbf{32}, AMS , Rhode Island, 1968


\bibitem{VK91} N.J. Vilenkin, A.U. Klimyk: Representations of Lie
groups and special functions, Vol 1. Kluwer Academi Publisher,
Dorderecht, 1991


\bibitem {W79}N. Wallach, \emph{
The analytic continuation of the discrete series, I, II}, Trans.
Amer. Math. Soc. \textbf{251} (1979), 1-17; 19-37.

\end{thebibliography}
\end{document}